\documentclass[11pt]{article}

\usepackage{latexsym,color,graphicx,amsmath,amssymb,amsthm}
\usepackage{url}

\newtheorem{theorem}{Theorem}
\newtheorem{lemma}[theorem]{Lemma}
\newtheorem{proposition}[theorem]{Proposition}
\newtheorem{corollary}[theorem]{Corollary}

\def\reals{{\mathbb R}}

\def\cplx{{\mathbb C}}
\def\C{{\mathcal C}}
\def\O{{\mathcal O}}

\def\P{{\mathbb P}}
\def\sph{{\mathbb S}}
\def\deg{{\mathsf{deg}}}

\newcommand{\Deg}{D} 
\makeatletter
\newcommand{\ProofEndBox}{{\ifhmode\unskip\nobreak\hfil\penalty50 \else
          \leavevmode\fi\quad\vadjust{}\nobreak\hfill$\Box$
            \finalhyphendemerits=0 \par}}
\makeatother
\newcommand{\fl}[1]{\mathsf{FL}_{#1}}

\newcommand{\proofend}{\ProofEndBox\smallskip}

\begin{document}

\title{Incidences between points and lines on a two-dimensional variety\thanks{%
Work on this paper by Noam Solomon and Micha Sharir was
supported by Grant 892/13 from the Israel Science Foundation.
Work by Micha Sharir was also supported
by Grant 2012/229 from the U.S.--Israel Binational Science Foundation,
by the Israeli Centers of Research Excellence (I-CORE)
program (Center No.~4/11), and
by the Hermann Minkowski-MINERVA Center for Geometry
at Tel Aviv University.
}}

\author{
Micha Sharir\thanks{%
School of Computer Science, Tel Aviv University,
Tel Aviv 69978, Israel.
{\sl michas@post.tau.ac.il} }
\and
Noam Solomon\thanks{%
School of Computer Science, Tel Aviv University,
Tel Aviv 69978, Israel.
{\sl noam.solom@gmail.com} }
}

\maketitle

\begin{abstract}
We present a direct and fairly simple proof of the following
incidence bound: Let $P$ be a set of $m$ points and $L$ a set of $n$
lines in $\reals^d$, for $d\ge 3$, which lie in a common algebraic two-dimensional
surface of degree $\Deg$ that does not contain any 2-flat, so that
no 2-flat contains more than $s \le \Deg$ lines of $L$. Then the
number of incidences between $P$ and $L$ is
$$
I(P,L)=O\left(m^{1/2}n^{1/2}\Deg^{1/2} +
m^{2/3}\min\{n,\Deg^{2}\}^{1/3}s^{1/3} + m + n\right).
$$
When $d=3$, this improves the bound of Guth and Katz~\cite{GK2} for this special case,
when $d$ is not too large.

A supplementary feature of this work is a review, with detailed proofs,
of several basic (and folklore) properties of ruled surfaces in three dimensions.
\end{abstract}

\section{Introduction}

Let $P$ be a set of $m$ distinct points and $L$ a set of $n$
distinct lines in $\reals^d$. Let $I(P,L)$ denote the number of
incidences between the points of $P$ and the lines of $L$; that is,
the number of pairs $(p,\ell)$ with $p\in P$, $\ell\in L$, and
$p\in\ell$. If all the points of $P$ and all the lines of $L$ lie in
a common 2-flat, then the classical Szemer\'edi--Trotter
theorem~\cite{SzT} yields the worst-case tight bound
\begin{equation} \label{inc2}
I(P,L) = O\left(m^{2/3}n^{2/3} + m + n \right) .
\end{equation}
This bound clearly also holds in $\reals^d$, for any $d$, by projecting the given
lines and points onto some generic 2-flat. Moreover, the bound will
continue to be worst-case tight by placing all the points and lines
in a common 2-flat, in a configuration that yields the planar lower
bound.

In the 2010 groundbreaking paper of Guth and Katz~\cite{GK2}, an improved
bound has been derived for $I(P,L)$, for a set $P$ of $m$ points and a set
$L$ of $n$ lines in $\reals^3$, provided that not too many lines of $L$ lie
in a common plane. Specifically, they showed:\footnote{%
  We skip over certain subtleties in their bound: They also assume that no
  \emph{regulus} contains more than $s=\sqrt{n}$ input lines, but then they are able
  also to bound the number of intersection points of the lines. Moreover,
  if one also assumes that each point is incident to at least three lines
  then the term $m$ in the bound can be dropped.}

\begin{theorem}[Guth and Katz~\cite{GK2}]
\label {ttt}
Let $P$ be a set of $m$ distinct points and $L$ a set of $n$ distinct lines
in $\reals^3$, and let $s\le n$ be a parameter,
such that
no plane contains more than $s$ lines of $L$. Then
$$
I(P,L) = O\left(m^{1/2}n^{3/4} + m^{2/3}n^{1/3}s^{1/3} + m + n\right).
$$
\end{theorem}
This bound (or, rather, an alternative formulation thereof) was a
major step in the derivation of the main result of \cite{GK2}, which
was an almost-linear lower bound on the number of distinct
distances determined by any finite set of points in the plane, a
classical problem posed by Erd{\H o}s in 1946~\cite{Er46}. Guth and Katz's
proof uses several nontrivial tools from algebraic and differential
geometry, most notably the Cayley--Salmon theorem on osculating
lines to algebraic surfaces in $\reals^3$, and various properties
of ruled surfaces. All this machinery comes on top of the major
innovation of Guth and Katz, the introduction of the
\emph{polynomial partitioning technique}.

In this study, in contrast with the aforementioned work of Guth and
Katz~\cite{GK2}, and with the follow-up studies of Guth~\cite{Gu14}
and of Sharir and Solomon~\cite{SS3d}, we do not need to explicitly use
the polynomial partitioning method, because we assume that the points
and lines all lie in a common surface of degree $\Deg$,
making the polynomial partitioning step superfluous. Concretely, our
main result is a simple and mostly elementary derivation of the
following result.

\begin{theorem}
\label {th:main} Let $P$ be a set of $m$ distinct points and $L$ a
set of $n$ distinct lines in $\reals^d$, and let $2\le s \le \Deg$ be two
integer parameters, so that all the points and lines lie in a common
two-dimensional algebraic variety $V$ of degree $\Deg$ that does not contain any
2-flat, and so that no 2-flat contains more than $s$ lines of $L$.
Then
\begin{equation} \label{st}
I(P,L) = O\left(m^{1/2}n^{1/2}\Deg^{1/2} + m^{2/3}\min\{n,\Deg^2\}^{1/3}s^{1/3} + m + n\right) .
\end{equation}
\end{theorem}

We assume in the theorem that $s$ is at most $\Deg$, but in fact
this assumption is superfluous and can be dropped. Indeed, for any
2-flat $\pi$, the intersection $\pi \cap V$ is a one-dimensional
plane algebraic curve of degree $\Deg$ in $\pi$ (this holds since
$V$ does not contain any 2-flat), and can therefore contain at most
$\Deg$ lines.

We also have the following easy and interesting corollary.
\begin{corollary}
\label {cor:main} Let $P$ be a set of $m$ distinct points and $L$ a
set of $n$ distinct lines in $\reals^d$, such that all the points and
lines lie in a common two-dimensional algebraic variety of constant
degree that does not contain any 2-flat. Then
$I(P,L) = O\left(m + n\right)$,
where the constant of proportionality depends on the degree of the surface.
\end{corollary}

For $d=3$, the corollary is not really new, because it can also be derived from
the analysis in Guth and Katz~\cite{GK2}, using a somewhat different
approach. Nevertheless, we are not aware of any previous explicit
statement of the corollary, even in three dimensions. As a matter of
fact, the corollary can also be extended (with a different bound though)
to the case where the containing
surface may have planar components. See a remark to that effect at the end
of the paper.

The significance of Theorem~\ref{th:main} is threefold:

\noindent
(a) In three dimensions, the bound improves the Guth--Katz bound when
$D\ll n^{1/2}$, for two-dimensional varieties $V$ that do not contain planes.
Note that the threshold $n^{1/2}$ is quite large---it is in fact larger than
the standard degree $O(m^{1/2}/n^{1/4})$ used in the analysis of
Guth and Katz~\cite{GK2} when $m\ll n^{3/2}$. We do not know how to
extend our results (or whether such an extension is at all possible)
to the case where $V$ does contain planes.

\noindent
(b) Another significant feature of our bound is that it does not contain
the term $nD$, which arises naturally in \cite{GK2} and other works,
and seems to be unavoidable
when $P$ is an arbitrary set of points. See an additional discussion of this
feature at the end of the paper.

\noindent
(c) It offers a sharp point-line incidence bound in arbitrary dimensions,
for the special case assumed in the theorem.

This paper is a step towards the study of point-line incidences on
arbitrary varieties, rather than just in $\reals^d$ (or $\cplx^d$).
Moreover, as has been our experience, the study of incidences in
$d$-space quickly reduces to questions about incidences within
lower-dimensional varieties, and tools for analyzing incidences of
this sort are very much in demand.

The analysis in this paper has indeed been used, as one of the key tools, in the
analysis in our companion paper~\cite{SS4d} on incidences between points and
lines in four dimensions. In this application, the lack of the term $nD$ is a crucial feature of
our result, which was required in the scenario considered in \cite{SS4d}.

The analysis in this paper makes extensive use of several properties of ruled surfaces
in $\reals^3$ or in $\cplx^3$. While these results exist as folklore in the literature,
we include here detailed and rigorous proofs thereof, making them more accessible
to the combinatorial geometry community. 



\section{Proof of Theorem~\ref{th:main}} \label{sec:pf1}

In most of the analysis in this section, we will consider
the case $d=3$. The reduction from an arbitrary dimension to $d=3$
will be presented at the end of the section.

\subsection{Preliminaries: Ruled surfaces} \label{sec:ruled}
The proof is based on several technical lemmas that establish various
properties of lines contained in an algebraic surface in three dimensions.
Some of them are also mentioned in our companion study~\cite{SS4d}
on point-line incidences in four dimensions, while others present
basic properties of ruled surfaces that are considered folklore
in the literature. Having failed to find rigorous proofs of
these properties, we provide here such proofs for the sake of completeness.

\paragraph{Singularity and flatness.}
The notion of singularities is a major concept, treated in full generality in
algebraic geometry (see, e.g., Kunz~\cite[Theorem VI.1.15]{Kunz} and
Cox et al.~\cite{CLO}), here we only recall some of their properties, and
only for a few special cases that are relevant to our analysis.

Let $V$ be a two-dimensional variety in $\reals^3$ or $\cplx^3$,
given as the zero set $Z(f)$ of some trivariate polynomial $f$.
Assuming $f$ to be square-free, a point $p\in Z(f)$ is
\emph{singular} if $\nabla f(p)=0$. For any point $p\in Z(f)$, let
$$
f(p+x)=f_{\mu}(x)+f_{\mu+1}(x)+\ldots
$$
be the Taylor expansion of $f$ near
$p$, where $f_j$ is the $j$-th order term in the expansion (which
is a homogeneous polynomial of $x$ degree $j$), and where we assume that
there are no terms of order (i.e., degree) smaller than $\mu$.
In general, we have $f_1(x) = \nabla f(p)\cdot x$,
$f_2(x) = \tfrac12 x^TH_f(p)x$, where $H_f$ is the Hessian matrix of $f$,
and the higher-order terms are similarly defined, albeit with more involved expressions.

If $p$ is singular,
we have $\mu\ge 2$. In this case, we say that $p$ is a singular point
of $Z(f)$ of \emph{multiplicity} $\mu$. For any point $p\in Z(f)$,
we call the hypersurface $Z(f_{\mu})$ the \emph{tangent cone} of
$Z(f)$ at $p$, and denote it by $C_p Z(f)$. If $\mu=1$, then $p$ is
non-singular and the tangent cone coincides with the (well-defned) tangent plane
$T_p Z(f)$ to $Z(f)$ at $p$.  We denote by $V_{sing}$ the locus of
singular points of $V$. This is a subvariety of dimension at most one;
see, e.g., Solymosi and Tao~\cite[Proposition~4.4]{SoTa}.

Similarly, let $\gamma$ be a one-dimensional algebraic curve in
$\reals^2$ or in $\cplx^2$, specified as $Z(f)$, for some bivariate
square-free polynomial $f$. Then $p\in Z(f)$ is \emph{singular} if
$\nabla f(p) = 0$. The \emph{multiplicity} $\mu$ of a point
$p\in \gamma$ is defined as in the three-dimensional case,
and we denote it as $\mu_{\gamma}(p)$; the multiplicity is at
least $2$ when $p$ is singular. The singular locus
$\gamma_{sing}$ of $\gamma$ is now a discrete set. Indeed, the fact that
$f$ is square-free guarantees that $f$ has no common factor with any
of its first-order derivatives, and B\'ezout's Theorem
(see, e.g., \cite[Theorem 8.7.7]{CLO}) then implies that
the common zero set of $f$, $f_x$, $f_y$, and $f_z$ is a (finite) discrete set.

By B\'ezout's Theorem, a line $\ell$ can intersect at most $D$
points of $\gamma$, \emph{counted with multiplicities}. To define this 
concept formally, as in, e.g., Beltrametti~\cite[Section 3.4]{BCGB},
let $\ell$ be a line and let $p\in\ell\cap\gamma$, such that $\ell$ 
is not contained in the tangent cone of $\gamma$ at $p$.
The \emph{intersection multiplicity} of $\gamma$ and $\ell$ at $p$ 
is the smallest order of a nonzero term of the Taylor expansion of
$f$ at $p$ in the direction of $\ell$.  The intersection multiplicity 
is also equal to $\mu_\gamma(p)$ (informally, this is the number of branches of
$\gamma$ that $\ell$ crosses at $p$, counted with multiplicity; 
see \cite[Section 8.7]{CLO} for a treatment on the intersection multiplicity in the plane).

We say that a line $\ell$ is a \emph{singular line} for $V$ if all of its
points are singular points of $V$. By Guth and Katz~\cite{GK} (see also
Elekes et al.~\cite[Corollary 2]{EKS}), the number of singular lines
contained in $V$ is at most $\deg(V)(\deg(V)-1)$.

Assume that $V=Z(f)$ is an irreducible algebraic surface. We say
that a non-singular point $x\in V$ is \emph{flat} if the second-order
Taylor expansion of $f$ at $x$ vanishes on the tangent plane
$T_x V$, or alternatively, if the \emph{second fundamental form} of $V$
vanishes at $x$ (see, e.g., Pressley~\cite{Pr}). Following Guth and
Katz~\cite{GK}, Elekes et al.~\cite[Proposition 6]{EKS} proved that
a non-singular point $x \in V$ is flat if and only if certain three
polynomials, each of degree at most $3\deg(V)-4$, vanish at $p$. A non-singular
line $\ell$ is said to be \emph{flat} if all of its non-singular points are flat.
By Guth and Katz~\cite{GK} (see also Elekes et al.~\cite[Proposition 7]{EKS}),
the number of flat lines fully contained in $V$ is at most $3\deg(V)^2-4\deg(V)$.

\paragraph{Ruled surfaces.}
For a modern approach to ruled surfaces, there are many references;
see, e.g., Hartshorne~\cite[Section V.2]{Hart83}, or
Beauville~\cite[Chapter III]{Beau}; see also Salmon~\cite{salmon}
and Edge~\cite{Edge} for earlier treatments of ruled surfaces.
We say that a real (resp., complex) surface $V$ is \emph{ruled by
real} (resp., \emph{complex}) \emph{lines} if every point $p$
in a Zariski-open\footnote{%
  The Zariski closure of a set $Y$ is the smallest algebraic variety
  $V$ that contains $Y$. $Y$ is Zariski closed if it is equal to its
  closure (and is therefore a variety), and is (relatively) Zariski open
  if its complement (within a given variety) is Zariski closed.
  See Cox et al.~\cite[Section 4.2]{CLO} for further details.}
dense subset of $V$ is incident to a real (complex) line that is fully
contained in $V$. This definition has been used in several
recent works, see, e.g.,~\cite{GK2,Kollar}; it is a slightly
weaker condition than the classical condition where it is required
that \emph{every} point of $V$ be incident to a line contained in
$V$ (e.g., as in~\cite{salmon}). Nevertheless, similarly to the proof of Lemma 3.4
in Guth and Katz~\cite{GK2}, a limit argument implies that the
two definitions are in fact equivalent. We give, in Lemma~\ref{below} below,
a short algebraic proof of this fact, for the sake of completeness.

We first recall the classical theorem of Cayley and Salmon.
Consider a polynomial $f \in \cplx[x,y,z]$ of degree $D\ge 3$. A
\emph{flecnode} of $f$ is a point $p\in Z(f)$ for which there exists a line
that passes through $p$ and \emph{osculates} to $Z(f)$ at $p$ to order three.
That is, if the direction of the line is $v$ then
$f(p) = 0$, and $\nabla_v f(p) = \nabla_v^2 f(p) = \nabla_v^3 f(p) = 0$,
where $\nabla_v f, \nabla^2_v f, \nabla^3_v f$ are, respectively, the
first, second, and third-order derivatives of $f$ in the direction $v$
(compare with the definition of singular points, as reviewed earlier,
for the explicit forms of $\nabla_v f$ and $\nabla^2_v f$).
The \emph{flecnode polynomial} of $f$, denoted $\fl{f}$, is the polynomial
obtained by eliminating $v$ from these three equations. As shown in
Salmon~\cite[Chapter XVII, Section III]{salmon}, the degree of $\fl{f}$ is
at most $11D-24$. By construction, the flecnode polynomial of $f$ vanishes on all
the flecnodes of $f$, and in particular on all the lines fully contained in $Z(f)$.

(Note that the correct formulation of Theorem~\ref{th:flec2} is over $\cplx$;
earlier applications, over $\reals$, as the one in Guth and Katz~\cite{GK2},
require some additional arguments to establish their validity; see 
Katz~\cite{Katz} for a discussion of this issue.)
\begin{theorem}[Cayley and Salmon~\cite{salmon}]
\label{th:flec2} Let $f \in \cplx[x,y,z]$ be a polynomial of degree
$\Deg\ge 3$. Then $Z(f)$ is ruled by (complex) lines if and only if
$Z(f) \subseteq Z(\fl{f})$.
\end{theorem}

\begin{lemma}
\label {le:rs} Let $f \in \cplx[x,y,z]$ be an irreducible
polynomial such that there exists a nonempty Zariski open dense set in $Z(f)$
so that each point in the set is incident to a line that is fully
contained in $Z(f)$. Then $\fl{f}$ vanishes identically on $Z(f)$,
and $Z(f)$ is ruled by lines.
\end {lemma}
\noindent{\bf Proof.}
Let $U\subset Z(f)$ be the set assumed in the lemma. By assumption
and definition, $\fl{f}$ vanishes on $U$, so $U$, and its Zariski closure,
are contained in $Z(f,\fl{f})$. Since $U$ is open, it must be two-dimensional.
Indeed, otherwise its complement would be a (nonempty)
two-dimensional subvariety of $Z(f)$
(a Zariski closed set is a variety). In this case, the complement must be equal to $Z(f)$,
since $f$ is irreducible, which is impossible since $U$ is nonempty.
Hence $Z(f,\fl{f})$ is also two-dimensional, and thus, by the same argument
just used, must be equal to $Z(f)$. Theorem~\ref{th:flec2} then implies that
$Z(f)$ is ruled by (complex) lines, as claimed. \proofend

\paragraph{Some additional tools from algebraic geometry.}
The main technical tool for the analysis is the following so-called
\emph{Theorem of the Fibers}.
Both Theorem~\ref{th:harr} and Theorem~\ref{th:harrdim}
hold (only) for the complex field $\cplx$.
\begin {theorem} [Harris~\protect{\cite[Corollary 11.13]{Har}}]
\label {th:harr} Let $X$ be a projective variety and $\pi: X \to \P^d$
be a homogeneous polynomial map (i.e., the coordinate functions
$x_0\circ \pi,\ldots,x_d \circ \pi$ are homogeneous polynomials); let
$Y=\pi(X)$ denote the image of $X$. For any $p\in Y$, let
$\lambda(p)=\dim(\pi^{-1}(\{p\}))$. Then $\lambda(p)$ is an upper
semi-continuous function of $p$ in the Zariski topology on $Y$; that
is, for any $m$, the locus of points $p\in Y$ such that
$\lambda(p)\ge m$ is Zariski closed in $Y$. Moreover, if $X_0 \subset X$ is
any irreducible component, $Y_0=\pi(X_0)$ its image, and $\lambda_0$
the minimum value of $\lambda(p)$ on $Y_0$, then
$$
\dim(X_0)=\dim(Y_0)+\lambda_0 .
$$
\end {theorem}

We also need the following theorem and lemma from Harris~\cite{Har}.
%

\begin {theorem} [Harris~\protect{\cite[Proposition 7.16]{Har}}]
\label{th:harrdim} Let $f: X \to Y$ be the map induced by the
standard projection map $\pi: \P^d \to \P^{r}$ (which retains $r$ of the
coordinates and discards the rest), where $r<d$, $X \subset \P^d$
and $Y\subset \P^{r}$ are projective varieties, $X$ is irreducible,
and $Y$ is the image of $X$. Then the general fiber\footnote{%
  The meaning of this statement is that the assertion holds for
  the fiber at any point outside some lower-dimensional
  exceptional subvariety.}
of the map $f$ is finite if and only if
$\dim(X)=\dim(Y)$. In this case, the number of points in a general
fiber of $f$ is constant.
\end {theorem}
In particular, when $Y$ is two-dimensional, there exist an integer
$c_f$ and an algebraic curve $\mathcal C_f \subset Y$, such that for
any $y \in Y\setminus \mathcal C_f$, we have $|f^{-1}(y)|=c_f$.
With the notations of Theorem~\ref{th:harrdim}, the set of points
$y\in Y$, such that the fiber of $f$ over $y$ is not equal to $c_f$
is a Zariski closed proper subvariety of $Y$. For more details, we refer
the reader to Shafarevich~\cite[Theorem II.6.4]{Shaf}, and to
Hartshorne~\cite[Exercise II.3.7]{Hart83}.

\begin {lemma} [Harris~\protect{\cite[Theorem 11.14]{Har}}]
\label {le:harr} Let $\pi: X\to Y$ be a polynomial map between two
projective varieties $X$, $Y$, with $Y=f(X)$ irreducible. Suppose that
all the fibers $\pi^{-1}(\{p\})$ of $\pi$, for $p\in Y$, are
irreducible and of the same dimension. Then $X$ is also irreducible.
\end {lemma}

\paragraph{Reguli.}
We rederive here the following (folklore) characterization of doubly
ruled surfaces. Recall that a \emph{regulus} is the surface spanned
by all lines that meet three pairwise skew lines.\footnote{%
  Technically, in some definitions (cf., e.g., Edge~\cite[Section I.22]{Edge})
  a regulus is a one-dimensional family of generator lines of the actual
  surface, i.e., a curve in the Pl\"ucker or Grassmannian space of lines,
  but we use here the notion of the surface spanned by these lines.}
\begin {lemma}
\label{doubly} Let $V$ be an irreducible ruled surface in $\reals^3$ or in $\cplx^3$
which is not a plane, and let $\C\subset V$ be an algebraic curve,
such that every non-singular point $p\in V\setminus \C$ is incident
to exactly two lines that are fully contained in $V$. Then $V$ is a regulus.
\end {lemma}

\begin{figure}[htb]
\begin{center}
\input{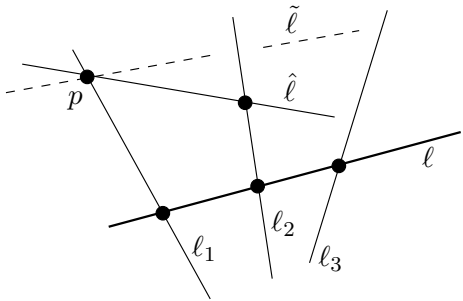}
\caption{The structure of $U_\ell$ in the proof of Lemma~\ref{doubly}.}
\label{freg}
\end{center}
\end{figure}

\noindent{\bf Proof.} As mentioned above (see also \cite{GK}), the
number of singular lines in $V$ is finite (it is smaller than
$\deg(V)^2$). For any non-singular line $\ell$, fully contained in
$V$, but not in $\C$, the union of lines $U_{\ell}$ intersecting
$\ell$ and fully contained in $V$ is a subvariety of $V$ (see Sharir
and Solomon~\cite[Lemma 8]{SS3d} 
for the easy proof). Each
non-singular point in $\ell \setminus \C$ is incident to another
line (other than $\ell$) fully contained in $V$, and thus $U_\ell$
is the union of infinitely many lines, and is therefore
two-dimensional. Since $V$ is irreducible, it follows that
$U_{\ell}=V$. Next, pick any triple of non-singular and
non-concurrent lines $\ell_1,\ell_2,\ell_3$ that are contained in
$V$ and intersect $\ell$ at distinct non-singular points of $\ell
\setminus \C$. There has to exist such a triple, for otherwise we
would have an infinite family of concurrent (or parallel) lines incident to $\ell$
and contained in $V$, and the plane that they span would then have
to be contained in $V$, contrary to assumption. See Figure~\ref{freg}
for an illustration. The argument given for $\ell$ applies equally well
to $\ell_1$, $\ell_2$, and $\ell_3$, and implies that 
$U_{\ell_1}=U_{\ell_2}=U_{\ell_3}=V$.
Assume that there exists some line $\tilde \ell\subset V$
intersecting $\ell_1$ at some non-singular point $p$, and that
$\tilde \ell \cap \ell_2= \emptyset$. (We treat lines here as
projective varieties, so this assumption means that $\tilde \ell$
and $\ell_2$ are skew to one another; parallel lines are considered
to be intersecting.) Since $p\in \ell_1 \subset V=U_{\ell_2}$, there
exists some line $\hat \ell$ intersecting $\ell_2$, such that
$\hat\ell \cap \ell_1 =\{p\}$. Hence there exist three distinct
lines, namely $\ell_1,\tilde \ell$ and $\hat \ell$, that are
incident to $p$ and fully contained in $V$. Since $p$ is
non-singular, it must be a flat point. Repeating this argument for
$3\deg(V)$ non-singular points $p \in \ell_1$, it follows that $\ell_1$ contains
at least $3\deg(V)$ flat points, and is therefore, by the properties
of flat points noted earlier, a flat line.
As is easily checked, $\ell_1$ can be taken to be an arbitrary
non-singular line among those incident to $\ell$, so it follows that every
non-singular point on $V$ is flat, and therefore, as shown in
\cite{EKS,GK}, $V$ is a plane, contrary to assumption.

Therefore, every non-singular line that intersects $\ell_1$ at a
non-singular point also intersects $\ell_2$, and, similarly, it
also intersects $\ell_3$. This implies that the intersection of
$V$ and the surface $R$ generated by the lines intersecting
$\ell_1,\ell_2$, and $\ell_3$ is two-dimensional, and is
therefore equal to $V$, since $V$ is irreducible.
Since $\ell_1,\ell_2$ and $\ell_3$ are pairwise skew,
$R=V$ is a regulus, as asserted.  \proofend

\paragraph{Real vs.~complex.}
At this point, we would like to elaborate about the field
over which the variety $V$ is defined. Most of the basic algebraic
geometry tools have been developed over the complex field $\cplx$,
and some care has to be exercised when applying them over the reals.
A major part of the theory developed in this paper is of this nature.
For example, both Theorems~\ref{th:harr} and~\ref{th:harrdim} hold
only over the complex field. As another important example,
one of the main tools at our disposal is the Cayley--Salmon
theorem (Theorem~\ref{th:flec2}), which applies over $\cplx$. Expanding
on a previously made comment, we note that even when $V$ is
a variety defined as the zero set of a real polynomial $f$, the
vanishing of the flecnode polynomial $\fl{f}$ only guarantees that
the set of complex points of $V$ is ruled by complex lines.

(A very simple example that illustrates this issue is the unit
sphere $\sigma$, given by $x^2+y^2+z^2=1$, which is certainly
not ruled by real lines, but the flecnode polynomial
of $f(x,y,z)=x^2+y^2+z^2-1$ (trivially) vanishes on $\sigma$.
This is the condition in the Cayley--Salmon theorem that
guarantees that $\sigma$ is ruled by (complex) lines, and indeed
it does, as is easily checked; in fact, for the same reason,
every quadric is ruled by complex lines.)

This issue has not been directly addressed in Guth and Katz~\cite{GK2},
although their theory can be adjusted to hold for the real case too,
as noted later in Katz~\cite{Katz}.

This is just one example of many similar issues that one must
watch out for.  It is a fairly standard practice in algebraic geometry
that handles a real algebraic variety $V$, defined by real polynomials,
by considering its complex counterpart $V_\cplx$, namely the set of
complex points at which the defining polynomials vanish. The rich
toolbox that complex algebraic geometry has developed allows one
to derive various properties of $V_\cplx$, but some care might be needed
when transporting these properties back to the real variety $V$,
as the preceding note concerning the Cayley--Salmon theorem illustrates.

In closing this discussion, we note that most of the results developed
in Section~\ref{ssec:pf} of this paper also apply over $\cplx$, except for one crucial step,
due to which we do not know how to extend Theorem~\ref{th:main} to the
complex domain. Nevertheless, we can derive a weaker variant of it for
the complex case---see a remark to that effect at the end of the paper.

\paragraph{Lines on a variety.}
In preparation for the key technical Theorem~\ref{singly}, given
below, we make the following comments. The \emph{Fano variety}
$F(V)$ of lines fully contained in a variety $V$ in three dimensions
is parametrized by the \emph{Pl\"ucker coordinates} of lines, as
follows (see, e.g., Griffiths and Harris~\cite[Section 1.5]{GrHa}).
For two points $x.y \in \P^3$, given in projective coordinates as
$x=(x_0,x_1,x_2,x_3)$ and $y=(y_0,y_1,y_2,y_3)$, let $\ell_{x,y}$
denote the (unique) line in $\P^3$ incident to both $x$ and $y$. The
Pl\"ucker coordinates of $\ell_{x,y}$ are given in projective
coordinates in $\P^5$ as
$(\pi_{0,1},\pi_{0,2},\pi_{0,3},\pi_{2,3},\pi_{3,1},\pi_{1,2})$,
where $\pi_{i,j}=x_iy_j-x_jy_i$. Under this parametrization, the set
of lines in $\P^3$ corresponds bijectively to the set of points in
$\P^5$ lying on the \emph{Klein quadric} given by the quadratic
equation $\pi_{0,1}\pi_{2,3}+\pi_{0,2}\pi_{3,1}+\pi_{0,3}\pi_{1,2} =
0$ (which is always satisfied by the Pl\"ucker coordinates of a
line). Given a surface $V$ in $\P^3$, the set of lines fully
contained in $V$, represented by their Pl\"ucker coordinates in
$\P^5$, is a subvariety of the Klein quadric, which, as above, is
denoted by $F(V)$, and is called the \emph{Fano variety} of $V$; see
Harris~\cite[Lecture 6, page 63]{Har} for details, and \cite[Example
6.19]{Har} for an illustration, and for a proof that $F(V)$ is
indeed a variety. The Pl\"ucker coordinates are continuous in the
sense that if one takes two points $\ell$, $\ell'$ on the Klein
quadric that are near each other, the lines in $\P^3$ that they
correspond to are also near to one another, in an obvious sense
whose precise details are omitted here. 

Given a plane $\pi$ by a homogeneous equation $A_0 x_0 + A_1 x_1 +
A_2 x_2 + A_3 x_3 = 0$, and a line $\ell$ not fully contained in
$\pi$, given in Pl\"ucker coordinates as
$(\pi_{0,1},\pi_{0,2},\pi_{0,3},\pi_{2,3},\pi_{3,1},\pi_{1,2})$, their point of
intersection is given in homogeneous coordinates by
$(A\cdot m, A \times m - A_0 d)$,
where $d=(\pi_{0,1},\pi_{0,2},\pi_{0,3})$, $m=(\pi_{2,3},\pi_{3,1},\pi_{1,2})$,
and where $\cdot$ stands for the scalar product, and $\times$ for
the vector product. This, together with the continuity argument
stated above, implies that, if the Fano variety $F(V)$ is
one-dimensional, and $\ell$ is a line represented by a regular point of $F(V)$,
then the cross section of the union of the lines that lie near $\ell$ in
$F(V)$ with a generic plane $\pi$ is a simple arc. When $\ell$ is a
singular point of $F(V)$, then the cross section of the union of the
lines that lie near $\ell$ in $F(V)$ with a generic plane $\pi$ is a
union of simple arcs meeting at $\ell \cap \pi$ where some of these 
arcs might appear with multiplicity; the number of these arcs is 
determined 
by the multiplicity of the singularity of $\ell$.

\paragraph{Singly ruled surfaces.}
The following result is another folklore result in the theory of
ruled surfaces, used in many studies (such as Guth and
Katz~\cite{GK2}). We give a detailed and rigorous proof, to make our
presentation as self-contained as possible; we are not aware of any 
similarly detailed argument in the literature.
\begin{theorem}
\label{singly} (a) Let $V$ be an irreducible ruled two-dimensional
surface of degree $\Deg>1$ in $\reals^3$ (or in $\cplx^3$), which is
not a regulus. Then, except for at most two exceptional lines, the
lines that are fully contained in $V$ are parametrized by an irreducible
algebraic curve $\Sigma_0$ in the Pl\"ucker space $\P^5$, and thus
yield a 1-parameter family of generator lines
$\ell(t)$, for $t\in \Sigma_0$, that depend continuously on the real or
complex parameter $t$. Moreover, if $t_1 \ne t_2$, and $\ell(t_1)
\ne \ell(t_2)$, then there exist sufficiently small and disjoint
neighborhoods $\Delta_1$ of $t_1$ and $\Delta_2$ of $t_2$, such that
all the lines $\ell(t)$, for $t\in \Delta_1\cup \Delta_2$, are
distinct.

\smallskip

\noindent
(b) There exists a one-dimensional curve $\C\subset V$, such that any point
$p$ in $V\setminus\C$ is incident to exactly one generator line of $V$.
\end{theorem}

\noindent{\bf Remark.} For a detailed description of the algebraic representation
of $V$ by generators, as in part (a) of the lemma, see Edge~\cite[Section II]{Edge}.

\medskip

\noindent{\bf Proof.} Assume first that we are working over $\cplx$.
Consider the Fano variety $F(V)$ of $V$, as defined above.
We claim that all the irreducible components of $F(V)$ are at most
one-dimensional. Informally, if any component $\Sigma_0$ of $F(V)$
were two-dimensional, then the set $\{ (p,\ell) \in V\times F(V)
\mid p\in\ell\}$ would be three-dimensional, so, ``on average'', the
set of lines of $F(V)$ incident to a point $p\in V$ would be
one-dimensional, implying that most points of $V$ are incident to
infinitely many lines that are fully contained in $V$, which can
happen only when $V$ is a plane (or a non-planar cone, which cannot
arise with a non-singular point $p$ as an apex), contrary to assumption.

To make this argument formal, consider the set (already mentioned above)
$$
W:=\{(p,\ell)\mid p \in \ell, \ell \in F(V)\}\subset V \times F(V) ,
$$
and the two projections
$$
\Psi_1 : W \to V, \quad \Psi_2: W \to F(V)
$$
to the first and second factors of the product $V\times F(V)$.

$W$ can formally be defined as the zero set of suitable homogeneous polynomials;
briefly, with an appropriate parameterization of lines in $\P^3$ and
the use of homogeneous coordinates, the
condition $p\in \ell$ can be expressed as the vanishing of two suitable
homogeneous polynomials, and the other defining polynomials are those that define
the projective variety $F(V)$. Therefore, $W$ is a projective variety.

Consider an irreducible component $\Sigma_0$ of $F(V)$ (which is also a
projective variety); put
$$
W_0:=\Psi_2^{-1}(\Sigma_0)=\{(p,\ell) \in W \mid \ell \in \Sigma_{0}\}.
$$
Since $W$ and $\Sigma_0$ are projective varieties, so is $W_0$.
As is easily verified, $\Psi_2(W_0)=\Sigma_0$ (that is, $\Psi_2$ is surjective).
We claim that $W_0$ is irreducible. Indeed, for any $\ell \in
\Sigma_0$, the fiber of the map $\Psi_2|_{W_0}: W_0 \to \Sigma_0$
over $\ell$ is $\{(p,\ell)\mid p \in \ell\}$ which is (isomorphic
to) a line, and is therefore irreducible of dimension one. As
$\Sigma_0$ is irreducible, Lemma~\ref{le:harr} implies that $W_0$
is also irreducible, as claimed.

For a point $p\in \Psi_1(W_0)$, consider the set
$\Sigma_{0,p}=\Psi_1|_{W_0}^{-1}(\{p\})$, put
$\lambda(p)=\dim(\Sigma_{0,p})$, and let
$\lambda_0 := \min_{p \in \Psi_1(W_0)} \lambda(p)$.
By the Theorem of the Fibers (Theorem~\ref{th:harr}), applied to the
map $\Psi_1|_{W_0}: W_0 \to V$, we have
\begin {equation}
\label {eq:fiber} \dim(W_0)= \dim(\Psi_1(W_0))+\lambda_0.
\end {equation}
We claim that $\lambda_0 = 0$. In fact, $\lambda(p) = 0$ for all
points $p\in V$, except for at most one point. Indeed, if $\lambda(p)\ge 1$
for some point $p\in V$, then $\Sigma_{0,p}$ is (at least) one-dimensional,
and $V$, being irreducible, is thus a cone with apex at $p$; since $V$ can
have at most one apex, the claim follows.  Hence $\lambda_0=0$, and therefore
\begin {equation} \label{eq:fibcon}
\dim(W_0) = \dim(\Psi_1(W_0)) \le \dim(V) = 2.
\end {equation}
Next, assume, for a contradiction, that $\dim(\Sigma_0)=2$. For a
point (i.e., a line in $\P^3$) $\ell \in \Psi_2(W_0)$, the set
$\Psi_2|_{W_0}^{-1}(\{\ell\})=\{(p,\ell)\mid p \in \ell\}$ is
one-dimensional (the equality follow from the way $W_0$ is defined). 
Conforming to the notations in the Theorem of the Fibers,
we have $\mu(\ell) := \dim\left(\Psi_2|_{W_0}^{-1}(\{\ell\}) \right) = 1$,
and thus $\mu_0:=\min_{\ell \in \Psi_2(W_0)} \mu(\ell)=1$.
Also, by assumption, $\dim(\Psi_2(W_0))=\dim(\Sigma_0)=2$.
By the Theorem of the Fibers, applied this time to
$\Psi_2|_{W_0}: W_0 \to \Sigma_0$, we thus have
\begin {equation}
\label {eq:fiber2}
\dim(W_0)= \dim(\Psi_2(W_0))+\mu_0 = 3,
\end {equation}
contradicting Equation~(\ref{eq:fibcon}). Therefore, every
irreducible component of $F(V)$ is at most one-dimensional, as claimed.

Let $\Sigma_0$ be such an irreducible component, and let
$W_0:=\Psi_2^{-1}(\Sigma_0)$, as above. As argued, for every
$p\in V$, the fiber of $\Psi_1|_{W_0}$ over $p$ is non-empty and
finite, except for at most one point $p$ (the apex of $V$ if $V$ is
a cone). Since $W_0$ is irreducible, Theorem~\ref{th:harrdim} implies
that there exists a Zariski open set $\mathcal O\subseteq V$,
such that for any point $p\in \mathcal O$,
the fiber of $\Psi_1|_{W_0}$ over $p$ has fixed cardinality $c_f$. Put
$\C:= V\setminus {\mathcal O}$. Being the complement of a Zariski open
subset of the two-dimensional irreducible variety $V$, ${\C}$ is (at
most) a one-dimensional variety. If $c_f \ge 2$, then, by
Lemma~\ref{doubly}, $V$ is a regulus. Otherwise, $c_f=1$ ($c_f$ cannot
be zero for a ruled surface), meaning that, for every $p\in
V\setminus {\C}$, there is exactly one line $\ell$, such that
$(p,\ell)\in W_0$, i.e., $\Sigma_0$ contains exactly one line
incident to $p$ and contained in $V$.

Moreover, we observe that the union of lines of $\Sigma_0$ is the entire variety $V$. 
Indeed, by Equation (\ref{eq:fiber}), it follows that $\dim(W_0)=\dim(\Psi_1(W_0))$. 
By Equations (\ref{eq:fibcon}) and (\ref{eq:fiber2}),
without the contrary assumption that $\dim(\Sigma_0)=2$, it follows that 
$\dim(W_0)=2$, and therefore $\dim(\Psi_1(W_0))=2$. That is, the variety
$\Psi_1(W_0)$, which is the union of the lines of $\Sigma_0$, must be the entire
variety $V$, because it is two-dimensional and is contained in the
irreducible variety $V$.


To recap, we have proved that if $\Sigma_0$ is a one-dimensional
component of $F(V)$, then the union of lines that belong to
$\Sigma_0$ covers $V$, and that there exists a one-dimensional
subvariety (a curve) $\C\subset V$ such that, for every
$p\in V\setminus {\C}$, $\Sigma_0$ contains exactly one line
incident to $p$ and contained in $V$.

Since $V$ is a ruled surface, some component of $F(V)$ has to be
one-dimensional, for otherwise we would only have a finite number of
lines fully contained in $V$. We claim that there is exactly one
irreducible component of $F(V)$ which is one-dimensional. Indeed,
assume to the contrary that $\Sigma_0,\Sigma_1$ are two (distinct)
one-dimensional irreducible components of $F(V)$. As we observed,
the union of lines parameterized by $\Sigma_0$ (resp.,
$\Sigma_1$) covers $V$. Let $\C_0$, $\C_1\subset V$ denote the
respective excluded curves, so that, for every
$p\in V\setminus \C_0$ (resp., $p\in V\setminus \C_1$) there
exists exactly one line in $\Sigma_0$ (resp., $\Sigma_1$) that
is incident to $p$ and contained in $V$.

Next, notice that the intersection $\Sigma_0 \cap \Sigma_1$ is a
subvariety strictly contained in the irreducible one-dimensional
variety $\Sigma_0$ (since $\Sigma_0$ and $\Sigma_1$ are two distinct
irreducible components of $F(V)$), so it must be zero-dimensional,
and thus finite.  Let $\C_{01}$ denote the union of the finitely many
lines in $\Sigma_0 \cap \Sigma_1$, and put $\C:= \C_0\cup\C_1\cup\C_{01}$.
For any point $p\in V\setminus \C$, there are two (distinct) lines incident
to $p$ and fully contained in $V$ (one belongs to $\Sigma_{0,p}$ and the
other to $\Sigma_{1,p}$). Lemma~\ref{doubly} (with $\C$ as defined above)
then implies that $V$ is a regulus, contrary to assumption.

In other words, the unique one-dimensional irreducible component
$\Sigma_0$ of $F(V)$ serves as the desired 1-parameter family of
generators for $V$. (The local parameterization of $\Sigma_0$ can be
obtained, e.g., by using a suitable Pl\"ucker coordinate of its lines.) 
In addition to $\Sigma_0$, there is a finite number of zero-dimensional
components (i.e., points) of $F(V)$. They correspond to a finite
number of lines, fully contained in $V$, and not parametrized by
$\Sigma_0$. Since the union of the lines in $\Sigma_0$ covers $V$,
any of these additional lines $\ell$ is \emph{exceptional}, in
the sense that each point on $\ell$ is also incident to a generator 
(different from $\ell$). By Guth and Katz~\cite[Corollary 3.6]{GK2}, 
there are at most two such exceptional lines, so there are at most two 
zero-dimensional components of $F(V)$. For the sake of completeness, we 
sketch a proof of our own of this fact; see Lemma~\ref{only2exc} in the appendix.

This establishes part (a) of the lemma, when $V$ is defined over $\cplx$.


If $V$ is defined over $\reals$, we proceed as above, i.e., consider
instead the complex variety $V_\cplx$ corresponding to $V$. As we
have just proven, the unique one-dimensional irreducible component
$\Sigma_0$ of $F(V)$ (regarded as a complex variety) is a (complex)
1-parameter family of generators for the set of complex points of
$V$. Since $V$ is real, the (real) Fano variety of $V$ consists of
the real points of $F(V)$, i.e., it is $F(V)\cap \P^5(\reals)$. As
we have argued, the (complex) $F(V)$ is the union of $\Sigma_0$ with
at most two other points. If $\Sigma|_\reals := \Sigma_0\cap\P^5(\reals)$ 
were zero-dimensional, the real $F(V)$ would also be discrete, so $V$
would fully contain only finitely many (real) lines, contradicting
the assumption that $V$ is ruled by real lines. Therefore, 
$\Sigma_0|_\reals$ is a one-dimensional irreducible component of the
real Fano variety of $V$. (It is irreducible, since otherwise the 
complex $\Sigma_0$ would be reducible too, as is easily checked.)



Summarizing, we have shown that there exists exactly one irreducible
one-dimensional component $\Sigma_0$ of $F(V)$, and a corresponding
one-dimensional subvariety ${\C}\subset V$, such that, for each
point $p\in V\setminus {\C}$, $\Sigma_0$ contains exactly one line
that is incident to $p$ (and contained in $V$). In addition to
$\Sigma_0$, $F(V)$ might also contain up to two zero-dimensional
(i.e., singleton) components, whose elements are the exceptional 
lines mentioned above. Let $\mathcal D$ denote the union of ${\C}$ 
and of the at most two exceptional lines; ${\mathcal D}$ is
clearly a one-dimensional subvariety of $V$. Then, for any point
$p \in V\setminus \mathcal D$, there is exactly one line incident to
$p$ and fully contained in $V$, as claimed. This establishes part
(b), and thus completes the proof of the lemma. \proofend


To complete the picture, we bridge between the ``conservative'' and
the ``liberal'' definitions of a ruled surface.
\begin {lemma} \label{below}
Let $V$ be a two-dimensional irreducible surface for which there exists a Zariski open
subset $\O\subseteq V$ with the property that each point $p\in\O$ is incident to
a line that is fully contained in $V$. Then this property holds for every point of $V$.
\end {lemma}

\noindent{\bf Proof.}
Consider the variety $W=\{(p,\ell)\mid p\in \ell, \ell \in F(V)\}$ used in the
proof of Theorem~\ref{salta}.
The projection $\Psi_1(W)$, in the notation of the preceding proof,
is a variety (assuming that we work in the complex projective setting)
contained in $V$. By assumption, it contains the Zariski open (and hence dense)
subset $\O$, so, being a variety (i.e., Zariski closed), it is the entire $V$, as
claimed.
\proofend

Ruled surfaces that are neither planes nor reguli, are called
\emph{singly ruled} surfaces (a terminology justified by 
Theorem~\ref{singly}). As in the proof of Theorem~\ref{singly}), a line $\ell$,
fully contained in an irreducible singly ruled surface $V$, such
that $\ell$ contains infinitely many points, each incident to
another line fully contained in $V$, is called an \emph{exceptional}
line of $V$. If there exists a point $p_V \in V$, which is incident
to infinitely many lines fully contained in $V$, then $p_V$ is
called an \emph{exceptional} point of $V$. By Guth and
Katz~\cite{GK2}, $V$ can contain at most one exceptional point $p_V$
(in which case $V$ is a cone with $p_V$ as its apex), and at most
two exceptional lines, as already noted. (This argument can be shown
to hold over both $\reals$ and $\cplx$.) 


\subsection{The proof}
\label{ssec:pf}

For a point $p$ on an irreducible singly ruled surface $V$, which is not the
exceptional point of $V$, we let $\Lambda_V(p)$ denote the number of
generator lines passing through $p$ and fully contained in $V$ (so if
$p$ is incident to an exceptional line, we do not count that line in
$\Lambda_V(p)$). We also put $\Lambda_V^*(p) := \max\{0,\Lambda_V(p)-1\}$.
Finally, if $V$ is a cone and $p_V$ is its
exceptional point (that is, apex), we put $\Lambda_V(p_V) = \Lambda_V^*(p_V):=0$.
We also consider a variant of this notation, where we are also given a
finite set $L$ of lines (where not all lines of $L$ are necessarily
contained in $V$), which does not contain any of the (at most two)
exceptional lines of $V$.  For a point $p\in V$, we let $\lambda_V(p;L)$
denote the number of lines in $L$ that pass through $p$ and are fully
contained in $V$, with the same provisions as above, namely that
we do not count incidences with exceptional lines, nor do we cound incidences
with an exceptional point, and put $\lambda_V^*(p;L) := \max\{0,\lambda_V(p;L)-1\}$.
If $V$ is a cone with apex $p_V$, we put $\lambda_V(p_V;L) = \lambda^*_V(p_V;L) = 0$.
We clearly have $\lambda_V(p;L) \le \Lambda_V(p)$ and
$\lambda^*_V(p;L) \le \Lambda^*_V(p)$, for each point $p$.


\begin{lemma}
\label{firstflip} Let $V$ be an irreducible singly ruled two-dimensional
surface of degree $\Deg>1$ in $\reals^3$ or in $\cplx^3$. Then, for
any line $\ell$, except for the (at most) two exceptional lines of
$V$, we have
\begin{align*}
& \sum_{p \in \ell \cap V} \Lambda_V(p) \le \Deg \quad\quad\text{if $\ell$ is not fully contained in $V$} , \\
& \sum_{p \in \ell \cap V} \Lambda^*_V(p) \le \Deg \quad\quad\text{if $\ell$ is fully contained in $V$} .
\end{align*}
\end{lemma}

\noindent{\bf Proof.}
We note that the difference between the two cases arises because we do not
want to count $\ell$ itself---the former sum would be infinite when $\ell$
is fully contained in $V$. Note also that if $V$ is a cone and $p_V\in\ell$,
we ignore in the sum the infinitely many lines incident to $p_V$ and
contained in $V$.

The proof is a variant of an observation due to
Salmon~\cite{salmon} and repeated in Guth and Katz~\cite{GK2} over
the real numbers, and later in Koll\'ar~\cite{Kollar} over general fields.

By Theorem~\ref{singly}(a), excluding the exceptional lines of $V$, the
set of lines fully contained in $V$ can be parameterized as a (real
or complex) 1-parameter family of generator lines $\ell(t)$,
represented by the irreducible curve $\Sigma_0 \subseteq F(V)$.
Let $V^{(2)}$ denote the locus of points of $V$ that are incident to at
least two generator lines fully contained in $V$. By
Theorem~\ref{singly}(b), $V^{(2)}$ is contained in some one-dimensional
curve $\C$.

Let $p\in V\cap\ell$ be a point incident to $k$ generator lines of $V$,
other than $\ell$, for some $k\ge 1$.
In case $V$ is a cone, we assume that $p \ne p_V$.  Denote the
generator lines incident to $p$ (other than $\ell$, if $\ell\subset
V$, in which case it is assumed to be a generator) as $\ell_i =
\ell(t_i)$, for $t_i\in\Sigma_0$ and for $i=1,\ldots,k$. 
(If $\ell_i$ is a singular point of $F(V)$,
it may arise as $\ell(t_i)$ for several values of $t_i$, and we pick 
one arbitrary such value.) Let $\pi$ be a generic
plane containing $\ell$, and consider the curve $\gamma_0 = V\cap
\pi$, which is a plane curve of degree $D$. Since
$V^{(2)}\subseteq\C$ is one-dimensional, a generic choice of $\pi$
will ensure that $V^{(2)}\cap\pi$ is a discrete set (since $\ell$ is
non-exceptional, it too meets $V^{(2)}$ in a discrete set).

There are two cases to consider: If $\ell$ is fully contained in $V$
(and is thus a generator), then $\gamma_0$ contains $\ell$. In this
case, let $\gamma$ denote $\gamma_0\setminus\ell$; it is also a
plane algebraic curve, of degree at most $\Deg-1$. Otherwise, we put
$\gamma := \gamma_0$. By Theorem~\ref{singly}(a), we can take, for
each $i=1,\ldots,k$, a sufficiently small open (real or complex)
neighborhood $\Delta_i$ along $\Sigma_0$ containing $t_i$, 
so that, for any $1 \le i < j \le k$, all the lines $\ell(t)$, 
for $t \in \Delta_i \cup \Delta_j$, are distinct. 
Put $V_i:= \bigcup_{t\in\Delta_i} \ell(t)$.
Recall that $V_i \cap \pi_i$ is either a simple arc
or a union of simple arcs meeting at $p$ (depending on whether or
not $\ell_i$ is a regular point of $\Sigma_0$); in the latter case,
take $\gamma_i$ to be any one of these arcs. Each of the arcs
$\gamma_i$ passes through $p$ and is contained in $\gamma$.
Moreover, since $\pi$ is generic, the arcs $\gamma_i$ are all
distinct. Indeed, for any $i\ne j$, and any point
$q\in\gamma_i\cap\gamma_j$, there exist $t_i \in \Delta_i, t_j \in
\Delta_j$ such that $\ell(t_i)\cap \pi= \ell(t_j)\cap \pi = q$, and
$\ell(t_i)\ne \ell(t_j)$ (by the properties of these neighborhoods).
Therefore, any point in $\gamma_i \cap \gamma_j$ is incident to (at
least) two distinct generator lines fully contained in $V$. Again,
the generic choice of $\pi$ ensures that
$\gamma_i\cap\gamma_j\subseteq V^{(2)}$ is a discrete set, so, in
particular, $\gamma_i$ and $\gamma_j$ are distinct.

We have therefore shown that (i) if $\ell$ is not contained in $V$
then $p$ is a singular point of $\gamma$ of multiplicity at least
$k$ (for $k\ge 2$; when $k=1$ the point does not have to be
singular), and (ii) if $\ell$ is contained in $V$ then $p$ is
singular of multiplicity at least $k+1$. We have $k \ge
\Lambda_V(p)$ (resp., $k \ge \Lambda^*_V(p)$) if $\ell$ is not fully
contained (resp., is fully contained) in $V$. (We may have an
inequality if $V$ is a cone and $p$ is its apex, since we then do
not count in $\Lambda_V(p)$ or in $\Lambda^*_V(p)$ the lines that
pass through $p$.) As argued at the beginning of Section~\ref{sec:ruled},
the line $\ell$ can intersect $\gamma$ in at most $\Deg$ points,
\emph{counted with multiplicity}, and the result follows.  \proofend

We also need the following result, established by Guth and
Katz~\cite{GK}; see also \cite{EKS}. It is an immediate consequence
of the Cayley--Salmon theorem (Theorem~\ref{th:flec2}) and a suitable
extension of B\'ezout's theorem for intersecting surfaces
(see Fulton~\cite[Proposition 2.3]{Fu84}).
\begin{proposition} [Guth and Katz \protect{\cite{GK2}}]
\label{caysal} Let $V$ be an irreducible two-dimensional variety in
$\cplx^3$ of degree $\Deg$. If $V$ fully contains more
than $11\Deg^2-24\Deg$ lines then $V$ is ruled by (possibly complex) lines.
\end{proposition}


\begin {corollary}
\label{co:caysal}
Let $V$ be an irreducible two-dimensional variety in
$\reals^3$ or $\cplx^3$ of degree $\Deg$ that does not contain any planes.
Then the number of lines that are fully contained in the union of the non-ruled
components of $V$ is $O(\Deg^2)$.
\end{corollary}

\noindent{\bf Proof.}
Let $V_1,\ldots,V_k$ denote those irreducible components of $V$ that are
not ruled by lines.  By Proposition~\ref{caysal}, for each $i$, the number
of lines fully contained in $V_i$ is at most $11\deg(V_i)^2-24\deg(V_i)$.
Summing over $i=1,\ldots,k$, the number of lines fully contained in the
union of the non-ruled components of $V$ is at most
$\sum_{i=1}^k 11\deg(V_i)^2= O(\Deg^2)$. \proofend

The following theorem, which we believe to be of independent interest in
itself, is the main technical ingredient of our analysis.
\begin{theorem} 
\label{salta} Let $V$ be a possibly reducible two-dimensional
algebraic surface of degree $D>1$ in $\reals^3$ or in $\cplx^3$,
with no linear components. Let $P$ be a set of $m$ distinct points on $V$ and
let $L$ be a set of $n$ distinct lines fully contained in $V$. Then there
exists a subset $L_0\subseteq L$ of at most $O(\Deg^2)$ lines, such
that the number of incidences between $P$ and $L\setminus L_0$
satisfies
\begin{equation} \label{lmlstar}
I(P,L\setminus L_0) = O\left(m^{1/2}n^{1/2}D^{1/2} + m + n\right) .
\end{equation}
\end{theorem}

\noindent{\bf Proof.} Consider the irreducible components
$W_1,\ldots,W_k$ of $V$. By Corollary~\ref{co:caysal}, the number of lines
contained in the union of the non-ruled components of $V$ is $O(\Deg^2)$,
and we place all these lines in the exceptional set $L_0$. In what follows
we thus consider only ruled components of $V$. For simplicity, continue
to denote them as $W_1,\ldots,W_k$, and note that $k\le D/2$.

We further augment $L_0$ as follows. We first dispose of lines of
$L$ that are fully contained in more than one ruled component $W_i$.
We claim that their number is $O(D^2)$. Indeed, for any pair $W_i$,
$W_j$ of distinct components, the intersection $W_i \cap W_j$ is a
curve of degree $\deg(W_i)\deg(W_j)$, which can therefore contain at
most $\deg(W_i)\deg(W_j)$ lines (by the generalized version of
B\'ezout's theorem~\cite[Proposition 2.3]{Fu84}, already mentioned in
connection with Proposition~\ref{caysal}).
Since $\sum_{i=1}^k \deg(W_i)\le D$, we have
$$
\sum_{i\ne j} \deg(W_i)\deg(W_j) \le \left(\sum_i \deg(W_i)\right)^2 = O(D^2) ,
$$
as claimed. We add to $L_0$ all the $O(D^2)$ lines in $L$ that are contained in
more than one ruled component, and all the exceptional lines of all
singly ruled components. The number of lines of the latter kind is at most
$2k \le 2\cdot(D/2) = D$, so the size of $|L_0|$ is still $O(D^2)$.
Hence, each line of $L_1:=L\setminus L_0$ is fully contained in a
unique ruled component of $V$, and is a generator of that component.

The strategy of the proof is to consider each line $\ell$ of $L_1$,
and to estimate the number of its incidences with the points of $P$ in
an indirect manner, via Lemma~\ref{firstflip}, applied to $\ell$ and
to each of the ruled components $W_j$ of $V$. We recall that $\ell$
is fully contained in a unique component $W_i$, and treat that component
in a somewhat different manner than the treatment of the other components.

In more detail, we proceed as follows.
We first ignore, for each singly ruled \emph{conic} component $W_i$, the incidences
between its apex (exceptional point) $p_{W_i}$ and the lines of $L_1$ that
are contained in $W_i$. We refer to these incidences as \emph{conical incidences}
and to the other incidences as \emph{non-conical}. When we talk about a line
$\ell$ incident to another line $\ell'$ at a point $p$, we will say that
$\ell$ is \emph{conically incident} to $\ell'$ (at $p$) if $p$ is the
apex of some conic component $W_i$ and $\ell'$ is fully contained in $W_i$
(and thus incident to $p$). In all other cases, we will say that $\ell$
is \emph{non-conically incident} to $\ell'$ (at $p$).
(Note that this definition is asymmetric in $\ell$ and $\ell'$; in
particular, $\ell$ does not have to lie in the cone $W_i$.)
We note that the number
of conical point-line incidences is at most $n$, because each line of $L_1$ is fully
contained in a unique component $W_i$, so it can be involved in at most
one conical incidence (at the apex of $W_i$, when $W_i$ is a cone).

We next prune away points $p\in P$ that are non-conically incident to at
most three lines of $L_1$. (Note that $p$ might be an apex of some conic
component(s) of $V$; in this case $p$ is removed if it is incident to
at most three lines of $L_1$ that are not contained in any of these
components.) We lose $O(m)$ (non-conical) incidences in this process.
Let $P_1$ denote the subset of the remaining points.

\noindent{\bf Claim.}
Each line $\ell\in L_1$ is non-conically incident, at points of $P_1$,
to at most $4\Deg$ other lines of $L_1$.

\noindent{\bf Proof.}
Fix a line $\ell\in L_1$ and let $W_i$ denote the unique ruled component
that fully contains $\ell$. Let $W_j$ be any of the other ruled components.
We estimate the number of lines of $L_1$ that are non-conically incident
to $\ell$ and are fully contained in $W_j$.

If $W_j$ is a regulus, there are at most four such lines, since $\ell$ meets
the quadratic surface $W_j$ in at most two points, each incident to exactly
two generators (and to no other lines contained in $W_j$).
Assume then that $W_j$ is singly ruled.  By Lemma~\ref{firstflip}, we have
$$
\sum_{p \in \ell \cap W_j} \lambda_{W_j}(p;L_1) \le
\sum_{p \in \ell \cap W_j} \Lambda_{W_j}(p) \le \deg(W_j) .
$$
Note that, by definition, the above sum counts only non-conical incidences
(and only with generators of $W_j$, but the exceptional lines of $W_j$ have
been removed from $L_1$ anyway). For $W_j$ a regulus, we write the bound $4$ as $\deg(W_j)+2$.

We sum this bound over all components $W_j\ne W_i$, including the reguli.
Denoting the number of reguli by $\rho$, which is at most $D/2$, we obtain a total of
$$
\sum_{j\ne i} \deg(W_j) + 2\rho \le \Deg + 2\rho \le 2\Deg .
$$
Consider next the component $W_i$ containing $\ell$. Assume first that $W_i$
is a regulus. Each point $p\in P_1\cap\ell$ can be incident to at most one other line
of $L_1$ contained in $W_i$ (the other generator of $W_i$ through $p$).
Since $p$ is in $P_1$, it is non-conically incident to at least $3-2=1$
other line of $L_1$, contained in some other ruled component of $V$.
That is, the number of lines that are (non-conically) incident to
$\ell$ and are contained in $W_i$, which apriorily can be arbitrarily large,
is nevertheless at most the number of other lines (not contained in $W_i$)
that are non-conically incident to $\ell$, which, as shown above, is at most $2\Deg$.

If $W_i$ is not a regulus, Lemma~\ref{firstflip} implies that
$$
\sum_{p \in \ell \cap W_i} \Lambda^*_{W_i}(p) \le \deg(W_i) \le \Deg ,
$$
where again only non-conical incidences are counted in this sum (and only
with generators).
That is, the number of lines of $L_1$ that are non-conically incident
to $\ell$ (at points of $P_1$) and are contained in $W_i$ is at most $\Deg$.
Adding the bound for $W_i$, which has just been shown to be either $D$ or $2D$,
to the bound $2\Deg$ for the other components, the claim follows.
\proofend

To proceed, choose a threshold parameter $\xi\ge 3$, to be determined
shortly. Each point $p\in P_1$ that is non-conically incident to at most
$\xi$ lines of $L$ contributes at most $\xi$ (non-conical) incidences,
for a total of at most $m\xi$ incidences. (Recall that the overall
number of conical incidences is at most $n$.) For the remaining
incidences, let $\ell$ be a line in $L_1$ that is incident to $t$
points of $P_1$, so that each such point $p$ is non-conically incident
to at least $\xi+1$ lines of $L_1$ (one of which is $\ell$).
By the preceding claim, we have $t\le 4\Deg/\xi$; summing this
over all $\ell \in L_1$, we obtain a total of at most $4n\Deg/\xi$
incidences. We can now bring back the removed points of $P\setminus P_1$,
since the non-conical incidences that they are involved in are counted
in the bound $m\xi$. That is, we have
$$
I(P,L_1) \le m\xi + n + \frac{4n\Deg}{\xi} .
$$
We now choose $\xi = (n\Deg/m)^{1/2}$. For this choice to make sense,
we want to have $\xi\ge 3$, which will be the case if $9m\le n\Deg$.
In this case we get the bound $O\left(m^{1/2}n^{1/2}\Deg^{1/2} + n\right)$.
If $9m>n\Deg$ we take $\xi=3$ and obtain the bound $O(m)$.
Combining these bounds, the lemma follows. \proofend

\paragraph{The final stretch.}
It remains to bound the number $I(P,L_0)$ of incidences involving
the lines in $L_0$. We have $|L_0|=O(\Deg^2)$; actually,
$|L_0| = \min\{n,O(\Deg^2)\}$. We estimate $I(P,L_0)$ using
Guth and Katz's bound (\cite{GK2}; see Theorem~\ref{ttt}), recalling
that no plane contains more than $s$ lines of $L_0$. We thus obtain
\begin{align} \label{lstar}
I(P,L_0) & =
O\left(m^{1/2}|L_0|^{3/4}+m^{2/3}|L_0|^{1/3}s^{1/3}+m+|L_0|\right) \\
& = O\left(m^{1/2}n^{1/2}\Deg^{1/2} +
m^{2/3}\min\{n,\Deg^2\}^{1/3}s^{1/3} + m + n \right) \nonumber .
\end{align}
Combining the bounds in (\ref{lmlstar}) and in (\ref{lstar}) yields the asserted
bound on $I(P,L)$.

\paragraph{Reduction to three dimensions.}
To complete the analysis, we need to consider the case where $V$
is a two-dimensional variety embedded in $\reals^d$, for $d>3$.

Let $H$ be a generic 3-flat, and denote by $P^*, L^*$, and $V^*$ the
respective projections of $P, L$, and $V$ onto $H$. Since $H$ is
generic, we may assume that all the projected points in $P^*$ are
distinct, and so are all the projected lines in $L^*$. Clearly,
every incidence between a point of $P$ and a line of $L$ corresponds
to an incidence between the projected point and line. Since no
2-flat contains more than $s$ lines of $L$, and $H$ is generic,
repeated applications of Theorem~\ref{ap:th} in the Appendix imply 
that no plane in $H$ contains more than $s$ lines of $L^*$.

One subtle point is that the set-theoretic projection $V^*$ of $V$
does not have to be a real algebraic variety (it is only a semi-algebraic
set), but it is always contained in a two-dimensional real algebraic variety
$\tilde V$, which we call, as we did in an earlier work~\cite{SSsocg}, the
\emph{algebraic projection} of $V$; it is the zero set of all polynomials
belonging to the ideal of polynomials vanishing on $V$, after eliminating
variables in the complementary space of $H$ (this is also known as an
\emph{elimination ideal} of $V$; see Cox et al.~\cite{CLO} for details).
We can also think of $\tilde V$
as the Zariski closure of $V^*$. Since the closure of a projection does not
increase the original degree (see, e.g., Harris~\cite{Har}), $\deg(\tilde V)\le D$.
That $\tilde V$ does not contain a 2-flat follows by a suitable adaptation
of the argument in Sharir and Solomon~\cite[Lemma 2.1]{SSsocg} 
(which is stated there for $d=4$ over the reals), that applies for general $d$ and
over the complex field too.

In conclusion, we have $I(P,L) \le I(P^*,L^*)$, where $P^*$ is a set of
$m$ points and $L^*$ is a set of $n$ lines, all contained in the
two-dimensional algebraic variety $\tilde V$, which is of degree at
most $D$ and does not contain any plane, and no plane contains more than
$s$ lines of $L^*$. The preceding analysis thus implies that the bound
asserted in the theorem applies in any dimension $d\ge 3$.
\proofend

\section{Discussion}

\noindent{\bf (1)}
We note that most of the proof of Theorem~\ref{th:main} can be carried out
over the complex domain. The only place where we (implicitly) use the fact
that the underlying field is $\reals$ is in the final step, where we apply the
bound of Guth and Katz, as a ``black box'', to $I(P,L_0)$. If we skip this step,
we remain with the weaker, albeit still useful, Theorem~\ref{salta},
which holds over $\cplx$ too. This is the part of our analysis that we apply
in our companion work~\cite{SS4d} on point-line incidences in four dimensions.

\medskip

\noindent{\bf (2)}
As mentioned in the introduction, Corollary~\ref{cor:main} can be extended to
the case where $V$, which is of constant degree $\Deg$, also contains planes.
Here too, we assume that no plane contains more than $s$ lines of $L$, but
this time it is not necessarily the case that $s\le \Deg$.

Let $\pi_1,\ldots,\pi_k$ denote the planar components of $V$, where $k\le\Deg=O(1)$.
For each $i=1,\ldots,k$, the number of incidences within $\pi_i$, namely, between the
set $P_i$ of points contained in $\pi_i$ and the set $L_i$ of lines fully contained
in $\pi_i$, is
$$
I(P_i,L_i) = O\left( |P_i|^{2/3}|L_i|^{2/3} + |P_i| + |L_i| \right)
= O\left( m^{2/3}s^{2/3} + m + s \right).
$$
Summing these bounds over the $k=O(1)$ planes, we get the same asymptotic bound for
the overall number of the incidences within these planes. Any other incidence
between a point $p$ lying in one of these planes $\pi_i$ and a line $\ell$
not contained in $\pi_i$ can be uniquely identified with the intersection of $\ell$
with $\pi_i$. The overall number of such intersections is at most $nk = O(n)$.
This leads to the following extension of Corollary~\ref{cor:main}.
\begin{corollary}
\label {cor:mainx} Let $P$ be a set of $m$ distinct points and $L$ a
set of $n$ distinct lines in $\reals^d$, and let $s\le n$ be a parameter,
such that all the points and lines lie in a common algebraic surface of
constant degree, and no 2-flat contains more than $s$ lines of $L$. Then
$$
I(P,L) = O\left(m^{2/3}s^{2/3} + m + n\right) ,
$$
where the constant of proportionality depends on the degree of the surface.
\end{corollary}

\medskip

\noindent{\bf (3)}
To us, one of the significant achievements of the analysis in this paper is
that our bound does not include the term $O(n\Deg)$. Such a term arises, as in the
preceding remark, when one considers incidences between points lying in some
irreducible component of $V$ and lines not contained in that component. These
incidences can be bounded by $n\Deg$, by charging them, as above, to
line-component intersections. This term also arises naturally in the
analysis of Guth and Katz~\cite{GK2} when considering incidences involving
lines not fully contained in the zero set of the partitioning polynomial
that is used in their analysis.
When $\Deg$ is large (say, $\Deg=O(n^{1/2})$), the term $n\Deg$ is too
large when compared with the bound of \cite{GK2}, but the term
$m^{1/2}n^{1/2}\Deg^{1/2}$ is fine.

\medskip

\noindent{\bf (4)}
Another issue that arises in our analysis is that, in the final step of the
proof, we bound the number of incidences involving lines of $L_0$ using the
Guth--Katz bound. It would be interesting, aesthetically speaking, to
replace this step by a different, simpler, and more direct derivation.
A more profound challenge would be to make this part work also over $\cplx$,
so that our bound would hold over the complex field too.

\medskip

\noindent{\bf (5)}
Concerning lower bounds, we do not know whether our result is tight.
A specific subproblem here is to bound the number of incidences within a
singly ruled irreducible surface (in three dimensions, say). The analysis
in Guth and Katz~\cite{GK2} yields the bound $O(m+nD)$, which (see item
(3) above) is in general too weak.

\medskip

\noindent{\bf (6)}
Finally, an interesting challenge is to establish a similar bound for $I(P,L)$,
for the case where the points of $P$ lie on a two-dimensional variety $V$,
but the lines need not be contained in $V$.

\appendix

\section{Some proofs}

\paragraph{Generic projections preserve non-planarity.}
\begin {theorem}
\label{ap:th} 
Let $\ell_1,\ell_2,\ell_3$ be three non-coplanar lines in $\reals^d$. 
Then, under a generic projection of $\reals^d$ onto some hyperplane $H$, 
the respective images $\ell^*_1,\ell^*_2,\ell^*_3$ of these lines
are still non-coplanar.
\end {theorem}

\noindent {\bf Proof.} 
Assume without loss of generality that the (generic) hyperplane $H$ 
onto which we project passes through $(0,0,0,0)\in \reals^d$, 
and let $w$ denote the unit vector normal to $H$. The projection 
$h:\;\reals^d\mapsto H$ is then given by $h(v) = v-(v\cdot w) w$.

Assume first that two of the three given lines, say $\ell_1,\ell_2$,
are skew (i.e., not coplanar). Let $\tilde \ell_1,\tilde \ell_2$ 
denote their projection onto $H$. If $\tilde \ell_1, \tilde \ell_2$ 
are coplanar they are either intersecting or parallel. If
they are intersecting, then there are points $p_1 \in \ell_1, p_2 \in \ell_2$ 
that project to the same point, i.e., $p_1-p_2$ has the
same direction as $w$. Then $w$ belongs to the set
$\{\frac {p_1-p_2} {\|p_1-p_2\|} \mid p_1 \in \ell_1, p_2 \in \ell_2\}$.
Since this is a two-dimensional set, it will be avoided for a generic choice of $w$,
which is a generic point in $\sph^{d-1}$, a set that is at least three-dimensional.

If $\tilde \ell_1,\tilde \ell_2$ are parallel, let $v_1,v_2$ denote the
directions of $\ell_1,\ell_2$. Since $v_1-(v_1\cdot w) w$ and
$v_2-(v_2\cdot w) w$ are vectors in the directions of
$\tilde \ell_1,\tilde \ell_2$, and are thus parallel, it follows that $w$ 
must be a linear combination of $v_1$ and $v_2$. Since $\|w\|=1$, the 
resulting set of possible directions is only one-dimensional, and, 
again, it will be avoided with a generic choice of $w$.

We may therefore assume that every pair of lines among $\ell_1,\ell_2,\ell_3$
are coplanar. Since these three lines are not all coplanar, the only two 
possibilities are that either they are all mutually parallel, or all concurrent. 

Assume first that they are concurrent, say they all pass through the origin.
Their projections are in the directions
$v_i-(v_i\cdot w) w$, for $i=1,2,3$. If these projections are coplanar then there
exist coefficients $\alpha_1,\alpha_2,\alpha_3$, not all zero, such that 
$\sum_i \alpha_i (v_i-(v_i\cdot w) w) = 0$. That is, putting
$u:= \sum_i \alpha_iv_i$, we have $u=(u\cdot w)w$, so $u$ is proportional
to $w$. In this case $w$ belongs to the set
$\{\frac{\sum_i \alpha_i v_i} {\|\sum_i \alpha_i v_i\|} \mid 
\alpha_1, \alpha_2, \alpha_3 \in \reals \text{ or } \cplx \}$. Again, being
a two-dimensional set, it will be avoided by a generic choice of $w$.

In the remainng case, the lines $\ell_1$, $\ell_2$, $\ell_3$ are mutually 
parallel, i.e., they all have the same direction $v$. Put, as above, 
for $i=1,2,3$, $\ell_i=\{p_i+tv\}_t$, and choose $p_i$ so that $p_i\cdot v = 0$.
The plane $\pi_0$ spanned by $p_1,p_2,p_3$ is projected to the plane $\pi$
spanned by the points $p_i^* = p_i - (p_i\cdot w)w$, for $i=1,2,3$ (since
$p_1$, $p_2$, $p_3$ are not collinear, they will not project into collinear 
points in a generic projection), and
the three lines project into a common plane if and only if their projections
are fully contained in $\pi$, meaning that the projection
$v^* = v-(v\cdot w)w$ is parallel to $\pi$, so it must be a linear 
combination of $p_1^*$, $p_2^*$, and $p_3^*$. A similar argument to those 
used above shows that a generic choice of $w$ will avoid the resulting 
two-dimensional set of forbidden directions.

This completes the proof. \proofend

Note that the ``hardest'' case in the preceding theorem is $d=4$. As $d$ 
increases, it becomes easier to avoid the forbidden two-dimensional sets of directions.

In the analysis in Section~\ref{ssec:pf}, the goal is to project $\reals^d$ onto
some generic 3-flat so that non-coplanar triples of lines do not project to
coplanar triples. This is easily achieved by repeated applications of Theorem~\ref{ap:th},
reducing the dimension one step at a time.

\paragraph{Exceptional lines on a singly ruled surface.}
In view of the proofs of Theorem~\ref{singly} and Lemma~\ref{below},
every point on such a surface $V$ is incident to at least one generator. Hence
an exceptional (non-generator) line is a line $\ell\subset V$ such that every
point on $\ell$ is incident to a generator (different from $\ell$).
\begin {lemma} \label{only2exc}
Let $V$ be an irreducible ruled surface in $\reals^3$ or in $\cplx^3$, which 
is neither a plane nor a regulus. Then $V$ contains at most two exceptional lines.
\end {lemma}

\noindent{\bf Proof.} 
We use the property, established in Sharir and Solomon~\cite{SS3d},
that for a line $\ell$ fully contained in $V$, the union $\tau(\ell)$ of 
the lines that meet $\ell$ and are fully contained in $V$ is a variety in
the complex projective space $\P^3(\cplx)$. Moreover, if $\ell$ is an 
exceptional line of $V$, then it follows by~\cite[Lemma 8]{SS3d} that
$\tau(\ell)=V$. (Indeed, $\tau(\ell)$ must be two-dimensional, since 
otherwise it would consist of only finitely many lines. Since $V$ is 
irreducible, $\tau(\ell)$ must then be equal to it.)

If $V$ contained three exceptional lines, $\ell_1, \ell_2$ and
$\ell_3$, then $V$ must be either a plane or a regulus. Indeed, otherwise, by
Theorem~\ref{singly} (whose proof does not depend on the number of 
exceptional lines), there exists a one-dimensional curve $\C\subset V$,
such that every point $p\in V\setminus \C$ is incident to
exactly one line $\ell_p$ fully contained in $V$.
As $p \in V\setminus\C$ and $\sigma(\ell_i)=V$, for $i=1,2,3$,
it follows that $\ell_p$ intersects $\ell_1,\ell_2$, and $\ell_3$.

If $\ell_1,\ell_2$, and $\ell_3$ are pairwise skew,
$p$ belongs to the regulus $R_{\ell_1,\ell_2,\ell_3}$ 
of all lines intersecting $\ell_1,\ell_2,$ and $\ell_3$. 
We have thus proved that $V\setminus \mathcal C$ is contained in
$R_{\ell_1,\ell_2,\ell_3}$, and as $R_{\ell_1,\ell_2,\ell_3}$
is irreducible, it follows that $V = R_{\ell_1,\ell_2,\ell_3}$. 

If $\ell_1,\ell_2$, and $\ell_3$ are concurrent but not coplanar then,
arguing similarly, $V$ is a cone with their common intersection point
as an apex. Since a (non-planar) cone has no exceptional lines, as is 
easily checked, we may ignore this case. 

Finally if any pair among $\ell_1,\ell_2$, $\ell_3$, say $\ell_1$, $\ell_2$,
are parallel then $V$ must be the plane that they span, contrary to
assumption. If $\ell_1$ and $\ell_2$ intersect at a point $\xi$, disjoint 
from $\ell_3$, then $V$ is the plane spanned by $\xi$ and $\ell_3$,
again a contradiction. 

Having exhausted all possible cases, the proof is complete. \proofend

\end{document}